\documentclass[12pt]{amsart}
\usepackage[mathscr]{eucal}
\usepackage{amsfonts}
\usepackage{amssymb,amsmath,amsthm}
\theoremstyle{plain}
\newtheorem{theorem}{Theorem}
\newtheorem{lemma}[theorem]{Lemma}

\theoremstyle{definition}

\begin{document}

\title[Infinitely Often Dense Bases]{Infinitely Often Dense Bases\\ of Integers with a Prescribed\\ Representation Function}
\author{Jaewoo Lee}
\address{Department of Mathematics\\ Borough of Manhattan Community College\\ The City University of New York\\ 199 Chambers Street\\ New York, NY 10007}
\email{jlee@bmcc.cuny.edu}
\date{February 8, 2007}
\begin{abstract}
 Nathanson constructed asymptotic bases for the integers with a prescribed representation function, then asked how dense they can be. We can easily obtain an upper bound using a simple argument. In this paper, we will see this is indeed the best bound we can get for asymptotic bases for the integers with an arbitrary representation function prescribed. 
\end{abstract}

\maketitle

\section{Introduction}\label{S:intro}

We will use the following notations: For sets $A, B$ of integers, any integer $t$ and a positive integer $h$, we define the \emph{sumset}
$ A + B = \{ a+b : a\in A, b \in B\}$, 
 the \emph{translation}
 $A + t = \{ a + t : a \in A\},$ and the \emph{dilation} 
$ h*A = \{ ha : a \in A\}$. 
Let $\mathbb{N}_0$ be the set of nonnegative integers. Then we define the representation function of $A$ as
\[r_A(n) = \mbox{card} \{(a,b) : a,b \in A,\; a \le b,\;
a+b=n \}, \] where $n$ is an integer. A set of integers A is called an \emph{additive basis} for integers if $r_A(n) \ge 1$ for all integers $n$. If all but finitely many  integers can be written as a sum of two elements from $A$, then $A$ is called an \emph{asymptotic basis} for integers. A set of integers
$A$ is called an \emph{unique representation basis} for integers if $r_A(n)=1$ for all integers $n$. Also, the \emph{counting function} for the set $A$ is \[ A(y,x) = \mbox{card} \{ a \in A : y \le a \le x \} \] for real numbers $x$ and $y$.

For bases of integers, Nathanson~\cite{NathSpa,NathInv} obtained the following:

\begin{theorem}[\cite{NathSpa, NathInv}]\label{T:Nath}
Let $ f \colon \mathbb{Z} \to \mathbb{N}_0 \cup \{\infty\}$ be any function such that the set $f^{-1}(0)$ is a finite set.
 \begin{enumerate}
 \item Let \, $\phi \, \colon \mathbb{N}_0 \to \mathbb{R}$ \; be any nonnegative function such that \; $\lim_{x \to \infty} \phi (x) = \infty$. Then there exist uncountably many asymptotic bases $A$ of integers such that $r_A(n) = f(n)$ for all integers \,$n$ and $A(-x,x) \le \phi(x)$ for all $x \ge 0$.
 \item There exist uncountably many asymptotic bases $A$ of integers such that $r_A(n) = f(n)$ for all integers \,$n$ and $A(-x,x) \gg x^{1/3}$ for all sufficiently large $x$.
 \end{enumerate}
\end{theorem}

Note that, under our assumption, a finite set $f^{-1}(0)$ just means $A$ is an asymptotic basis. Cilleruelo and Nathanson~\cite{Per} later improved the exponent ${1/3}$ in the second statement of  Theorem~\ref{T:Nath} to $\sqrt{2} -1 +o(1)$. Also,  \L uczak and Schoen~\cite{Luc} proved the second statement of Theorem~\ref{T:Nath} by showing that a set with a condition of Sidon type can be extended to a unique representation bases. 

When the representation function is arbitrarily given, the first statement of Theorem~\ref{T:Nath} means an asymptotic basis for the integers can be as sparse as we want, and the second statement of Theorem~\ref{T:Nath} means we can achieve a certain thickness. Therefore, we want to consider the following general question: \emph{Given an arbitrary representation function, what are the possible thickness for asymptotic bases for the integers?} This question was posed by Nathanson~\cite{NathInv}. 

We can obtain an upper bound quite easily, as shown in \cite{Uniq}. Let $A$ be any set of integers with a bounded representation function $r_A(n) \le r$ for some $r > 0$ for all integers $n$. Take $k = A(-x,x)$. Then there are $\frac{k(k+1)}{2}$ ways to make $a_i + a_j$\,, where $a_i, a_j$ are in $A$ with $|a_i|, |a_j| \le x$. All these sums belong to the interval $[-2x, 2x]$ and each number in that interval is represented by $a_i + a_j$ at most $r$ times. Therefore, \[ \frac{k(k+1)}{2} \le r(4x+1) \] and solving this for $k = A(-x,x)$ gives us 
$A(-x,x) \ll x^{1/2}$
 for all $x > 0$.

Since we are considering possible thickness for asymptotic bases with an arbitrary representation function given, the above argument provides us with an upper bound $x^{1/2}$ for the possible thickness. Then the next question we might ask is, can we find a better upper bound, or is $x^{1/2}$ the best upper bound we can get? One way to find a better upper bound is to find a representation function which require the thickness to have a better upper bound. 

We might want to start with unique representation bases. Nathanson~\cite{Uniq} asked the following: Does there exist a number $\theta < 1/2$ such that $A(-x,x) \le x^{\theta}$ for every unique representation basis $A$ and for all sufficiently large $x$? This question was answered negatively by Chen~\cite{Chen}.

\begin{theorem}[\cite{Chen}]\label{T:Chen}
For any \,$\epsilon > 0$, there exists an unique representation basis $A$ for the integers such that for infinitely many positive integers $x$, we have 
\[ A(-x,x) \ge x^{{1/2} - \epsilon} . \]
\end{theorem}

In this paper, we show that it is impossible to find a better upper bound for any representation function. Therefore, if we consider $A(-x,x)$ for infinitely many integers $x$ instead of all large $x$, $x^{1/2}$ is indeed the best upper bound. 
 \\

\textbf{Acknowledgment}. The author thanks Mel Nathanson for his very helpful advices.

\section{Preliminary Lemmas}\label{S:Lem}

From now on, $f$ will denote a function $f \colon \mathbb{Z} \to \mathbb{N}_0 \cup \{ \infty \}$ such that the set $f^{-1}(0)$ is a finite set. Then there exists a positive integer $d_0$ such that $f(n) \ge 1$ for all integers $n$ with $|n| \ge d_0$. Nathanson~\cite{NathSpa} proved the following:

\begin{lemma}[\cite{NathSpa}]\label{L:uk}
Given a function $f$ as above, there exists a sequence $U=\{ u_k \}_{k=1}^{\infty}$ of integers such that, for every $n \in \mathbb{Z}$ and $k \in \mathbb{N}$, 
\begin{equation}\label{E:uk}
f(n) = \mbox{card} \{ k \ge 1 : u_k = n \}.
\end{equation}
\end{lemma}

For the proof of Lemma~\ref{L:uk}, see Nathanson~\cite{NathSpa}.

\begin{lemma}\label{L:first}
Let $A$ be a finite set of integers with $r_A(n) \le f(n)$ for all integers \,$n$, \, $0 \notin A$, and for all integers $n$,\[r_A(n) \ge \# \{ i \le m: u_i = n\}\] for some integer \,$m$ which depends only on the set $A$. Then, there exists a finite set of integers $B$ such that $A \subseteq B$, $r_B(n) \le f(n)$ for all integers  $n$, \[r_B(n) \ge \# \{ i \!\le\! m+1 : u_i = n \}\] for all integers $n$, and  $0 \notin B$.
\end{lemma}

\begin{proof}
If $r_A(n) \ge \# \{ i \le m+1 : u_i = n\}$ for all $n$, then take $B=A$ and we are done. Otherwise, note that \[\#\{ i \le m : u_i = n\} = \# \{ i \le m+1 : u_i = n \}\] for all $n \ne u_{m+1}$ and \[\#\{ i \le m : u_i = u_{m+1}\} + 1 = \# \{ i\le m+1 : u_i = u_{m+1}\}.\] If \[r_A(n) < \#\{ i\le m+1 : u_i = n\}\] for some $n$ as we are assuming now, since \[r_A(n) \ge \# \{ i \le m : u_i = n\}\] for all $n$, we must have \[r_A(u_{m+1}) < \#\{ i \le m+1 : u_i=u_{m+1} \} \le f(u_{m+1}).\]

Let $d = \mbox{max}\{ d_0,\, |u_{m+1}|,\,|a|\ \mbox{where }a \in A \}$. Choose $c > 4d$ if $u_{m+1} \ge 0$ and $c < -4d$ if $u_{m+1} < 0$. Note that $|c| >4d$.

Let $B = A \cup \{ -c,\ c+u_{m+1}\}$. Then $2B$ has three parts:
\[ 2A, \ A+\{-c,\ c+u_{m+1}\}, \ \{ -2c,\ u_{m+1},\ 2c+2u_{m+1} \}. \]
If $a \in 2A$, then $-2d \le a \le 2d$. If $a \in A$, then if $u_{m+1} \ge 0$,\ we have $c > 0$, thus
\begin{eqnarray*}
a-c & \le & d-4d = -3d, \\ 
a + c + u_{m+1} & \ge & -d+4d+u_{m+1} \ge 3d +u_{m+1} \ge 3d\,,
\end{eqnarray*}
and if $u_{m+1} <0$,\ we have $c<0$, thus
\begin{eqnarray*}
a-c &\ge& -d+4d = 3d, \\
a+c+u_{m+1} &\le& d-4d+u_{m+1}=-3d+u_{m+1} \le -3d\,.
\end{eqnarray*}
Therefore, \[ 2A \,\cap\, A\! +\! \{-c,\ c+u_{m+1} \} = \emptyset\,, \] each element of $A\! +\! \{-c,\ c+u_{m+1}\}$ has an unique representation in the form of $a\! +\! \{-c, \ c+u_{m+1}\},\ a\in A$, and the same is true for $\{-2c,\ u_{m+1},\ 2c+2u_{m+1}\}$ \,in the form of $\{-c,\ c+u_{m+1}\} + \{-c,\ c+u_{m+1}\}$. And
\begin{eqnarray*}
|-2c| &=& 2|c| > 8d, \\
|2c+2u_{m+1}|&=& 2|c+u_{m+1}| \ge 2|c| \ge 8d\,,
\end{eqnarray*}
thus $2A \,\cap\, \{-2c, 2c+2u_{m+1}\} = \emptyset$ \,(recall that $c$ and $u_{m+1}$ has the same sign).
Also note that $A+\{-c,\ c+u_{m+1}\}$ and $\{-2c,\ u_{m+1},\ 2c+2u_{m+1}\}$ are disjoint. To see this, for example, if $a-c=2c+2u_{m+1}$ for some $a \in A$, then $a=3c+2u_{m+1}$ \,so $|a|=|3c+2u_{m+1}| \ge |3c| > 12d$, giving us a contradiction. Other cases are similar. 

Thus, we have
\[ r_B(n) = 
\left\{\begin{array}{cl}
 r_A(n) + 1 & \mbox{if } n=u_{m+1} \\
 r_A(n)     & \mbox{if } n\in \,2A \setminus \{u_{m+1}\} \\
 1          & \mbox{if } n\in \,2B \setminus \bigl\{2A \cup \{u_{m+1}\}\bigr\}.
\end{array}\right. \]

Now, we have $r_A(u_{m+1}) < f(u_{m+1})$, so\, \[r_B(u_{m+1}) = r_A(u_{m+1}) +1 \;\le f(u_{m+1}).\] And, if $n\in 2B \setminus \bigl\{2A \cup \{u_{m+1}\}\bigr\}$, then $|n| \ge d_0$, so $f(n) \ge 1 = r_B(n)$. Thus, $r_B(n) \le f(n)$ for all $n$.

Now, \,$r_A(u_{m+1}) \ge \#\{i\le m : u_i = u_{m+1}\}$ so 
\begin{eqnarray*}
r_A(u_{m+1})+1 & \ge & \#\{i\le m : u_i=u_{m+1}\} +1 \\
 & = & \#\{i \le m+1 : u_i = u_{m+1} \}, \end{eqnarray*}
therefore, 
\[ r_B(u_{m+1}) \ge \#\{ i \le m+1 : u_i = u_{m+1} \}. \]
If $n \in 2A \setminus \{u_{m+1}\}$, then \[r_B(n) = r_A(n) \ge \#\{i \le m : u_i = n \} = \#\{ i \le m+1 : u_i = n\}.\]
If $n \in 2B \setminus \bigl\{2A \cup \{u_{m+1}\}\bigr\}$, then 
\[ 0=r_A(n) \ge \#\{ i \le m : u_i = n\} = \#\{ i \le m+1 : u_i = n\}\]
so
\[ 0=\#\{ i \le m+1 : u_i = n \} \le 1 = r_B(n). \]
Thus, \[r_B(n) \ge \#\{ i \le m+1 : u_i =n\} \] for all $n$.
\end{proof}

\begin{lemma}\label{L:second}
Let $A$ be a finite set of integers with $r_A(n) \le f(n)$ for all $n$, and $0 \notin A$. Let $\phi (x) \colon \mathbb{N}_0  \to \mathbb{R}$ be a nonnegative function such that $\lim_{x \to \infty} \phi (x) = \infty$. Then for any $M > 0$, there exists an integer $x > M$ and a finite set of integers $B$ with $0 \notin B$, $A \subseteq B$, $r_B(n) \le f(n)$ for all $n$, and $B(-x,x) > \sqrt{x}/\phi(x)$.
\end{lemma}

\begin{proof}
It is well-known that there exists a Sidon set $D \subseteq [1,n]$ such that $|D| = n^{1/2} + o(n^{1/2})$ for all $n \ge 1$. Choose an integer $x$ which satisfies the following:
\begin{enumerate}
 \item $\phi (x) > M+ \sqrt{20T}$ where $T= \mbox{max} \{d_0,\ |a|\ \mbox{where }a \in A\},$\label{ff}
 \item $x$ is a multiple of $5T$,\label{fs}
 \item $x > M$,\label{ft}
 \item If $n$ is large enough, $|D|= \sqrt{n} +o(\sqrt{n}) > \sqrt{n}/2$. Let $x$ be large enough so that $n=x/5T$ would be large enough to  satisfy the above.\label{ffo}
\end{enumerate}

Let $B = A \cup \{ 5Td : d \in D\}$ where $D \subseteq [1,n]$ with $n= x/5T$ as above. Then 
\begin{eqnarray*}
B(-x,x) &\ge& B(0,x)\ge D(1,n)=|D| \\
  &> &\frac{\sqrt{n}}{2} = \frac{1}{2}\sqrt{\frac{x}{5T}} = \frac{\sqrt{x}}{\sqrt{20T}} > \frac{\sqrt{x}}{\phi(x) - M} > \frac{\sqrt{x}}{\phi (x)}. \end{eqnarray*}

Now, note that $2B$ has three parts:
\[ 2A, \ A+5Td \ \mbox{for } d\in D,\ \mbox{and } 5T(d_1+d_2)\ \mbox{for } d_
1,d_2 \in D.\] 
As earlier, we have $2A \,\cap\, (A + 5T\negmedspace*\negmedspace D) = \emptyset$\; and $2A \,\cap\, \{ 5T(d_1+d_2) : d_1, d_2 \in D\} = \emptyset$.
Now, if $a+5Td_1 = 5T(d_2 + d_3)$ for $a \in A$, $d_i \in D$, then $|a| = 5T\,|d_2+d_3-d_1|$. If $|\,d_2+d_3-d_1|=0$, then $a=0 \in A$, a contradiction. And if $|\,d_2+d_3-d_1| \ge 1$, then $|\,a| \ge 5T$, a contradiction. Thus 
\[ (A+5T\negmedspace*\!D) \cap \{5T(d_1+d_2) : d_1, d_2 \in D\} = \emptyset.\]
If $a_1+5Td_1 = a_2 +5Td_2$ \,for $a_1, a_2 \in\! A$,\, $d_1, d_2 \in\! D$, then $|\,a_1 - a_2|=5T|\,d_1-d_2|$. As above, this cannot happen unless $a_1=a_2$ \,and $d_1 =d_2$. And if \,$5T(d_1+d_2)=5T(d_3+d_4)$ for $d_i \in D$ with $d_1 \le d_2, \, d_3 \le d_4$, then $d_1+d_2=d_3+d_4$. Since $D$ is a Sidon set, we have $d_1=d_3,\ d_2=d_4$. Thus we have 
\[ r_B(n) = 
\left\{\begin{array}{cl}
 r_A(n) & \mbox{if } n \in 2A \\
 1    & \mbox{if } n\in 2B \setminus 2A.
\end{array}\right. \]

If $n \in 2B \setminus 2A$, \, $n \ge T \ge d_0$, so $f(n) \ge 1 = r_B(n)$.  Thus, $r_B(n) \le f(n)$ for all $n$.
\end{proof}

\section{Main Result}\label{S:Main}

\begin{theorem}\label{T:Main}
Let $f \colon \mathbb{Z} \to \mathbb{N}_0 \cup \{ \infty\}$ be a function with a finite set $f^{-1}(0)$. Let $\phi \colon \mathbb{N}_0 \to \mathbb{R}$ be any nonnegative function with $\lim_{x \to \infty} \phi(x) = \infty$. Then, there exists an asymptotic basis $A$ of integers such that $r_A(n) = f(n)$ for all $n$, and for infinitely many positive integers\, $x$, we have $A(-x,x) > \sqrt{x}/\phi(x)$.
\end{theorem}

\begin{proof}
Recall that given such a function $f$, we have $d_0$ and $\{ u_k\}$ as defined in Section~\ref{S:Lem}. 
We use induction to get an infinite sequence of finite sets of integers $A_1 \subseteq A_2 \subseteq \cdots$ and a sequence of positive integers $\{ x_i\}_{i=1}^{\infty}$ with $x_{i+1} > x_i$\, such that, for all positive integers $l$, we have 
\begin{enumerate}
 \item $r_{A_l}(n) \le f(n)$ for all $n$,\label{indfi}
 \item $r_{A_{2l}}(n),\, r_{A_{2l+1}}(n) \;\ge\, \#\{ i\le l\!+\!1 : u_i=n\}$ for all $n$,\label{indse}
 \item  $A_{2l-1}(-x_l, x_l) > \sqrt{x_l}/\phi(x_l)$,\label{indth}
 \item $0 \notin A_l$.\label{indfo}
\end{enumerate}

If $u_1 \ge 0$, take $c=4d_0 > 0$. If $u_1 <0$, take $c=-4d_0 <0$. Let $\alpha =|2c+2u_1| > 0$\;(\,thus $\alpha > |c|\,,|u_1|\,,d_0$). As before, if $n$ is large enough, there exists a Sidon set $D \subseteq [1,n]$ such that $|D| > \sqrt{n}/2$. Take such an integer $n$ which also satisfies $\phi(3\alpha n) > 2\!\sqrt{3\alpha}$.

Take $A_1 = 3\alpha\!*\!D \;\cup\, \{-c,\ c+u_1\}$ \,and \,$x_1=3\alpha n$. Then $2A_1$ has three parts:
\[ 2(3\alpha\!*\!D),\  \ 3\alpha\!*\!D+\{-c, \ c+u_1\},\  \ \{-2c, \ u_1,\ 2c+2u_1\}. \]

It's easy to see that these are pairwise disjoint(for example, if $3\alpha d -c=2c+2u_1$\,, then $3\alpha d=3c+2u_1$\,, so $3\alpha \le |\,3c+2u_1| < |\,4c+4u_1|=2\alpha$\,, \,a contradiction). If $3\alpha d_1 -c=3\alpha d_2+c+u_1$, then $2c+u_1=3\alpha (d_1-d_2)$. If $d_1 \ne d_2$, then $|2c+u_1| \ge 3\alpha$ but $|\,2c+u_1| < |\,2c+2u_1| = \alpha$. So $d_1\!=\!d_2$\,. Then $-c=c+u_1$\,, \,so $-2c=u_1$\,, a contradiction. Thus,  
\[ r_{A_1}(n) = 
\left\{\begin{array}{cl}
 1  & \mbox{if } n \in 2A_1 \\
 0  & \mbox{if } n \notin 2A_1.
\end{array}\right. \]

Now, $3\alpha d-c\ge 3\alpha - \alpha=2\alpha$\,, and also\, $3\alpha d+c+u_1 \ge 3\alpha -\alpha-\alpha = \alpha$\,. So if $n \!\in\! 2A_1 \!\setminus\! \{u_1\}$ then $|n| \ge d_0$\,, so $f(n) \ge 1$. And if $n= u_1$\,, then by the definition of $\{u_k\}$, we have $f(u_1)= \#\{ k : u_k=u_1 \,\} \ge 1$. Thus, for all $n \in 2A_1$\,,\  $r_{A_1}(n) =1 \le f(n)$. If $ n \notin 2A_1$\,, \, $r_{A_1}(n)=0 \le f(n)$. Therefore, for all $n$, $r_{A_1}(n) \le f(n)$. Now, we have $1=r_{A_1}(u_1) \ge \#\{ i \le 1 : u_i = u_1\}$. For other $n \ne u_1$\,, \ $r_{A_1}(n) \ge \#\{ i \le 1 : u_i = n\}=0$. Thus, $r_{A_1}(n) \ge \#\{i\le 1 : u_i=n\}$ for all $n$. Also, 
\begin{eqnarray*}
A_1(-x_1,\, x_1) &\ge&  A_1(1,\, 3\alpha n) \ge \,D(1,n) = |D| \\
               &> &   \frac{\sqrt{n}}{2}=  \frac{\sqrt{3\alpha n}}{2\!\sqrt{3\alpha}} > \frac{\sqrt{3\alpha n}}{\phi(3\alpha n)} = \frac{\sqrt{x_1}}{\phi(x_1)}\,.
\end{eqnarray*}
Thus $A_1$ satisfies all of the conditions~\eqref{indfi} to \eqref{indfo}.

Now, suppose we have $A_1 \subseteq A_2 \subseteq \cdots \subseteq A_{2l-1}$ and $x_1 < x_2 < \cdots < x_l$. By Lemma~\ref{L:first}, there exists $A_{2l}$ such that $A_{2l-1} \subseteq A_{2l}$ with $r_{A_{2l}}(n)\le f(n)$ for all $n$\, and 
\[ r_{A_{2l}}(n)\, \ge\, \#\{ i \le l\!+\!1 :\, u_i=n\} \] for all $n$\,, and $0 \notin A_{2\,l}$. Now, by Lemma~\ref{L:second}, there exists an integer $x_{l+\!1} > x_l$ and $A_{2l+\!1}$ with $0 \notin A_{2l+\!1}$\,, \ $A_{2l} \subseteq A_{2l+\!1}$\,, \ $r_{A_{2l+\!1}}(n) \le f(n)$ for all $n$\,, and \[A_{2l+\!1}(-x_{l+\!1}\,,\,  x_{l+\!1})\, > \,\frac{\sqrt{x_{l+\!1}}}{\phi(x_{l+\!1})} \,.\]
Also,\, $r_{A_{2l+\!1}}(n) \,\ge\, r_{A_{2l}}(n) \,\ge\, \#\{i \le l\!+\!1 : u_i=n\}$ \,for all $n$.

Now, let $A= \cup_{l=1}^{\infty} \;A_l$\,. By conditions \eqref{indfi} and \eqref{indse}\,, $r_A(n) = f(n)$ for all $n$ and 
\[  A(-x_k\,,\,x_k) \ge A_{2k-\!1}(-x_k\,,\,x_k) > \frac{\sqrt{x_k}}{\phi(x_k)} \]
for all $k$\,.
\end{proof}


\begin{thebibliography}{99}

\bibitem{Chen}
Y. Chen, \emph{A problem on unique representation bases}, European J. Combin. \textbf{28} (2007), 33-35
 \\
\bibitem{Per}
J. Cilleruelo and M. B. Nathanson, \emph{Dense sets of integers with prescribed representation functions}, preprint, 2006.
 \\
\bibitem{Luc}
T. \L uczak and T. Schoen, \emph{A note on unique representation bases for the integers}, Funct. Approx. Comment. Math. \textbf{32} (2004), 67-70.
 \\
\bibitem{Uniq}
M. B. Nathanson. \emph{Unique representation bases for the integers}, Acta Arith. \textbf{108} (2003), no. 1, 1-8.
 \\
\bibitem{NathInv}
M. B. Nathanson, \emph{The inverse problem for representation functions of additive bases}, in : Number Theory (New York, 2003), 253-262, Springer, New York, 2004.
 \\
\bibitem{NathSpa}
M. B. Nathanson, \emph{Every function is the representation function of an additive basis for the integers}, Port. Math. (N.S.) \textbf{62} (2005), no. 1, 55-72.

\end{thebibliography}
\end{document}